\documentclass[11pt,a4paper]{article}

\usepackage[utf8]{inputenc}
\usepackage[T1]{fontenc}
\usepackage{amsmath,amssymb,amsthm}
\usepackage{booktabs}
\usepackage{rotating}
\usepackage{graphicx}
\usepackage[margin=2.5cm]{geometry}
\usepackage{hyperref}
\usepackage{color}
\usepackage{pgfplots}
\pgfplotsset{compat=1.18}
\usepackage[english=british]{csquotes}

\usepackage{authblk}
\usepackage{orcidlink}

\newtheorem{theorem}{Theorem}
\newtheorem{lemma}[theorem]{Lemma}

\theoremstyle{definition}
\newtheorem{definition}[theorem]{Definition}
\theoremstyle{remark}

\newtheorem{example}[theorem]{Example}
\AtEndEnvironment{example}{\hfill$\diamond$}

\newcommand{\numNec}{m} \newcommand{\numSuf}{r} \newcommand{\Wmin}{\mathcal{W}^{\min}}
\newcommand{\Lmax}{\mathcal{L}^{\max}}
\newcommand{\LmaxX}{\mathcal{L}^{\max}_{\geq \numNec}}
\newcommand{\pop}{p}
\newcommand{\set}[1]{\left\{#1\right\}}
\newcommand{\V}{\mathcal{V}}
\newcommand{\W}{\mathcal{W}}
\newcommand{\Wz}{\mathcal{W}_{\numSuf}}
\newcommand{\LargeCountries}{B}
\newcommand{\SubsetLarge}{I}
\newcommand{\SmallCountries}{{\overline{B}}}
\newcommand{\SubsetSmall}{J}
\newcommand{\minPopWin}{\beta} \newcommand{\maxPopLose}{\alpha} \newcommand{\subsetsLargeCtoCover}{\mathcal{N}}
\newcommand{\subsetsLargeCtoCoverDom}{\mathcal{N}^*}
\newcommand{\subsetsLargeWithFinitAlpha}{\mathcal{U}}
\newcommand{\subsetsLargeWithFinitAlphaDom}{\mathcal{U}^*}
\newcommand{\bnb}{branch-and-bound }

\newcommand{\Z}{\mathbb{Z}}

\newcommand{\fitsize}[1]{\resizebox{\dimexpr\linewidth}{!}{#1}}

\newcommand{\ie}{i.e.,\ }

\title{Determining the dimension of the EU Council voting rule} \author{Eva~Ley~\orcidlink{0009-0000-7927-4757}}
\author{Maximilian~Merkert~\orcidlink{0000-0002-7838-445X}}
\date{} 

\affil{
TU~Braunschweig\\
Institute for Mathematical Optimization\\
Universitätsplatz 2, 38106 Braunschweig, Germany\\
\url{eva.ley@tu-braunschweig.de},~
\url{m.merkert@tu-braunschweig.de}
}

\begin{document}
\maketitle

\begin{abstract}
The dimension of a simple game is the smallest number of weighted games whose intersection it is. 
While most real-world voting rules have a relatively small dimension that is easy to determine, for the qualified majority rule of the Council of the European Union under the Treaty of Lisbon, with $28$ member states and the 2014 population data, the best published lower and upper bounds on the dimension are $8$ and $24$, respectively. 
Determining the exact dimension has been posed as a challenge to the community by Kurz and Napel in~2015.
We describe a heuristic for finding small systems of weighted games for simple games of the structure of the voting rule in the Treaty of Lisbon. 
Using our heuristic, we obtain an explicit representation of the rule as an intersection of eight weighted games with integer weights for the 2014 population data. 
This upper bound matches the previously best known lower bound from Kober and Weltge (2021). 
Thus, the dimension of the voting rule of the Council of the European Union under the Treaty of Lisbon for the 2014 population data is $8$.
\end{abstract}

\section{Introduction}

A \emph{simple game} on a finite player set $N$ is a family $\W \subseteq 2^N$ of \emph{winning coalitions} that is closed under
taking supersets; coalitions not in $\mathcal{W}$ are called \emph{losing}. 
It is \emph{weighted} if there are weights $v\in\mathbb{R}^N$ 
and a quota $t$ such that $S\in\mathcal{W}$ if and only if $v(S):=\sum_{i\in S}v_i\ge t$.
In short, we write $\W=[t;v]$ in this case.
Every simple game is an intersection of weighted games, and its \emph{dimension} is the smallest number of weighted games needed~\cite{TaylorZwicker1993}. The question of how many linear inequalities are required to separate any given subset of vertices of the $0$-$1$ hypercube from its remaining vertices was studied earlier in~\cite{Jeroslow75}, and corresponds to the concept of relaxation complexity, see~\cite{KW15}, restricted to binary points.
From~\cite{Jeroslow75} we have a general upper bound of $2^{|N|-1}$ on the dimension of any simple game. 
The dimension 
can indeed be exponential in the number of players~\cite{ODS17_GrowthDimCompleteSimpleGames}. The question \enquote{do there exist $k$ weighted games such that for all coalitions \dots} for $k$ polynomially bounded in the input size is in~$\Sigma_2^p$, the class typical of bilevel problems with an integer lower level, see~\cite{Woeginger21}.
Computing the dimension of a game given as an intersection of weighted games is NP-hard~\cite{DeinekoWoeginger2006}.
Verifying the equivalence of two weighted games is coNP-complete~\cite{EGGW08_DimensionalityVotingGames}.

Most voting rules used in practice are weighted, \ie of dimension one.
This might even hold if they are classically formulated as the intersection of several weighted games.
For example, in the United Nations Security Council a coalition wins if at least nine of its fifteen members vote in favour and none of the five permanent members vote against.
Ignoring abstentions, this can be formulated via a single weighted game with unit weights for all non-permanent members and a weight of 10 for all permanent members, combined with a quota of $54$, see~\cite{Woeginger21}, where the problem of determining the dimension of simple games and this very example are used as introduction to illustrate problems in $\Sigma_2^p$.
Examples of dimension two are the US federal legislative system and the procedure for amending the Canadian constitution~\cite{TaylorZwicker1993}.
The EU Council under the Treaty of Nice and the Legislative Council of Hong Kong have dimension three~\cite{Freixas02, CheungNg14}.

Since 1 November 2014, the Council of the EU decides according to the Treaty of Lisbon, applying a \emph{qualified majority} rule as the standard voting procedure, where a coalition $S \subseteq N$ wins if and only if
\begin{itemize}
    \item at least $55\%$ of the countries vote in favour and these countries represent at least $65\%$ of the population or
    \item at most three countries do not support a vote.
\end{itemize}
Thus, 
\begin{equation}\label{eq:EUrule-preBrexit}
  \bigl(|S|\ge 0.55 |N| \ \text{ and } \ \pop(S)\ge 0.65\,\pop(N)\bigr)
  \quad\text{or}\quad |S|\ge |N|-3
\end{equation}
where $N$ is the set of member countries and $\pop_i$ the population of country $i$.
For $|N|=28$, the rule is $\W_{EU28}=({\cal C}_{16}\cap {\cal P})\cup {\cal C}_{25}$ with the weighted games ${\cal C}_k=[k;1,\dots,1]$ and ${\cal P}=[0.65\,\pop(N);\pop]$.
Kurz and Napel~\cite{KurzNapel2016} showed $7\le\dim \W_{EU28}\le 13\,368$ for the $2014$ population data and posed determining the dimension as a challenge to the community. 
Chen, Cheung and Ng~\cite{ChenCheungNg2019} reduced the upper bound on the dimension to $24$ by giving a structural construction which represents ${\cal P\cup C}_{25}$ as an intersection of $23$ weighted games. 
Kober and Weltge~\cite{KoberWeltge2021} improved the lower bound to $8$ by exhibiting a set of $15$ losing coalitions whose associated hypergraph has a chromatic number of~$8$. 
Vossen~\cite{Vossen} tried to improve this lower bound using ideas from \cite{AHS23} but only reached $7$ as a lower bound.

For the $2014$ population data as in \cite{KurzNapel2016,ChenCheungNg2019}, reproduced in Table~\ref{tab:pop2014}, our main result is:
\begin{theorem}\label{thm:main}
For the $2014$ population data, the EU Council voting rule~\eqref{eq:EUrule-preBrexit} has dimension $8$.
\end{theorem}
As a certificate, we provide integer weight vectors and quotas for $8$ weighted games in Table~\ref{tab:games2014}. Note that verifying the equivalence of two weighted games is still a coNP-complete problem~\cite{EGGW08_DimensionalityVotingGames}. 
However, in our case the certificate can be verified in reasonable time even via enumeration of all minimal winning and maximal losing coalitions, using exact integer computations. We provide such a verification code via \cite{VerificationCode}.
Combined with the lower bound by Kober and Weltge~\cite{KoberWeltge2021}, the dimension of eight follows.
To the best of our knowledge, the lower bound of eight -- established to be tight here -- is the largest known dimension for a voting rule that is actually used by a major political institution. 

\paragraph{Outline}
We develop a heuristic to find systems of weighted games that represent games of the form $\W=({\cal C}_{\numNec} \cap {\cal P})\cup {\cal C}_{\numSuf}$ for integers $\numNec<\numSuf$, where $\cal P$ is a single weighted game. After introducing some preliminaries in Section~\ref{sec:preliminaries}, we describe our heuristic in Section~\ref{sec:heuristics}. For reducing the search space, we first introduce three assumptions on a potential system of weighted games.
These are heuristic in nature, and, in fact, all may result in not being able to find a system of dimension-many weighted games in general, as is demonstrated by small examples.
However, for the EU 2014 population data, it turns out  that one is able to close the gap under the heuristic assumptions.
First, fixing the weighted game ${\cal C}_{\numNec}$ to be included reduces the number of maximal losing coalitions to be considered.
Furthermore, we distinguish countries with larger and smaller population. Introducing the notion of $\LargeCountries$-reduced games, we find a system of weighted games for the EU Council where the weights of smaller countries are equal to their population.
In combination with the previous restrictions, assuming integer weights is a third restriction.
Finally, we optimally colour the resulting hypergraph with a \bnb approach.
This approach allows computing a system of dimension eight for the EU population data of 2014 in less than ten seconds, see Section~\ref{sec:results}.  
\section{Preliminaries}
\label{sec:preliminaries}

For a finite set of countries (players) $N$, let $p\in\mathbb{Z}^N_{\geq 0}$ (populations), $q\in\mathbb{Z}$, and integers $\numNec<\numSuf$ be given. 
We consider games of the form
\begin{equation}\label{eq:rule}
  S \in \W \iff \bigl(|S|\ge \numNec \text{ and } p(S)\ge q\bigr)\ \text{ or }\ |S|\ge \numSuf .
\end{equation}
We can write $\W$ as $\W={\cal C}_\numNec\cap({\cal P}\cup {\cal C}_\numSuf)$ with ${\cal C}_k=[k;1,\dots,1]$ and ${\cal P}=[q;p]$.
Note that this is a simple game, \ie every subset of a losing coalition is again losing and every superset of a winning coalition is winning.
Thus, it suffices to only consider the minimal winning and maximal losing coalitions explicitly.
Let $\Wmin$ denote the set of minimal winning and $\Lmax$ the set of maximal losing coalitions of the game~$\W$.
Furthermore, we can restrict the weights for weighted games to be nonnegative.
Let ${\cal V}_1,\dots,{\cal V}_k$ be weighted games with nonnegative weights.
We have $\W=\bigcap_{i=1}^k{\cal V}_i$ if and only if every $W\in\Wmin$ is winning in every ${\cal V}_i$, and every $L\in\Lmax$ is losing in at least one ${\cal V}_i$.

We order the countries by decreasing population, $\pop_1 \geq \pop_2 \geq \dots \geq \pop_{|N|}$. 
As the voting rule only depends on the size of the population and cardinalities of coalitions, this gives a simple characterization of coalitions in $\Wmin$ and $\Lmax$.
A winning coalition $W$ is minimal winning if and only if $W\setminus \{ \max(i \mid i \in W) \}$ is losing. 
A losing coalition $L$ is maximal losing if and only if $L\cup\{\min (j \mid j \notin L) \}$ is winning.

Consistent with \cite[Def.1]{KoberWeltge2021}, we call a set $\mathcal{S}\subseteq\Lmax$ \emph{separable} if a single weighted game~$\V$ accepts all of~$\Wmin$ and rejects all of $\mathcal{S}$, \ie $\V \supseteq \Wmin$ and $\V \cap \mathcal{S} = \emptyset$.
It is easy to see that separability is closed under taking subsets, and the dimension of $\W$ is the smallest number of separable sets covering $\Lmax$, \ie the chromatic number of the hypergraph on $\Lmax$ whose edges are the non-separable sets. 
Separability can be decided based on feasibility of a system of linear inequalities:

\begin{lemma}\label{lem:lp}
$\mathcal{S}\subseteq\Lmax$ is separable if and only if there are
$v\in\mathbb{R}^N_{\ge 0}$ and $t\in\mathbb{R}$ with
\begin{equation}\label{eq:lp}
  v(W)\ge t \ \ (W\in\Wmin),\qquad v(L)\le t-1 \ \ (L\in\mathcal{S}).
\end{equation}
In that case $[t-\tfrac12;v]$ is a weighted game with the required property.
\end{lemma}

\begin{proof}
A weighted game $[q;v]$ with $v(W)\ge q>v(L)$ for all $W\in\Wmin$, $L\in\mathcal{S}$ can be scaled so that the gap $\min_W v(W)-\max_L v(L)$ is at least $1$.
Set $t=\min_W v(W)$.
\end{proof}

\section{A Heuristic for finding upper bounds} \label{sec:heuristics}

\subsection{Fix one weighted game}

We fix ${\cal C}_{\numSuf}$ to be one of the weighted games in our representation.
This is also done in the construction for the best upper bound from the literature~\cite{ChenCheungNg2019}.
Then, only the coalitions in $\LmaxX:=\{L\in\Lmax:\ |L|\ge \numNec\}$ have to be covered by the remaining games.
Kober and Weltge mention that the losing coalition consisting of the largest $15=\numNec-1$ countries (by population) for EU 2014 data is particularly well-connected in the hypergraph they use for their proof of the lower bound of 8: This losing coalition forms a two-element non-separable set (corresponding to an edge in the hypergraph) with any losing coalition that consists of 23 or 24 countries~\cite{KoberWeltge2021}.

Fixing a game upfront can prevent us in general from reaching the minimal number of weighted games, as the following example shows.

\begin{example}\label{example:fixCxCanCost}
    Consider the game
    \begin{equation*}
        (|S| \geq 3 \text{ and } p(S) \geq 9) \text{ or } |S| \geq 5 \text{ with } p = (5, 3, 3, 1, 1, 1, 1),
    \end{equation*}
    which is equivalent to 
    $(1 \in S \text{ and } S \cap \set{2, 3} \neq \emptyset \text{ and } |S| \geq 3) \text{ or }  |S| \geq 5$.
    The maximal losing coalitions are $\Lmax = \set{ \set{1, 2}, \set{1, 3} } \cup \set{S : |S| = 4, 1 \notin S}$.
    We claim that it has a dimension of two which cannot be obtained if the weighted game ${\cal C}_3 = [3; 1, \dots, 1]$ is included.
    
    First, the dimension is at most two as the following two weighted games certify:
    \begin{align*}
        {\cal V}_1 &= [5;\,3,1,1,1,1,1,1]\text{, equivalent to }(1 \in S \text{ and } |S| \geq 3)\text{ or } |S| \geq 5\\
        {\cal V}_2 &= [5;\,1,5,5,1,1,1,1]\text{, equivalent to }(S \cap \set{2,3} \neq \emptyset) \text{ or }|S| \geq 5.
    \end{align*}
    Their intersection is exactly the described game: 
    ${\cal V}_1$ rejects $\set{1,2}$, $\set{1,3}$ and all losing coalitions of size $4$ without player $1$ and ${\cal V}_2$ rejects the four coalitions that contain $1$ but neither $2$ nor $3$.

    However, any system containing ${\cal C}_3$ needs at least three games: 
    ${\cal C}_3$ rejects $\set{1,2}$ and $\set{1,3}$ while all maximal losing coalitions of size $4$ must be rejected by the other games.
    This includes in particular the losing coalitions $L_1=\set{2,3,4,5}$ and $L_2=\set{1, 4, 6, 7}$, which are, however, not separable: 
    Consider the winning coalitions $W_1=\set{1,2,4}$ and $W_2=\set{3, 4, 5, 6, 7}$. 
    For all weights $v$, we have 
    \begin{equation*}
        v(L_1)+v(L_2) = v(W_1)+v(W_2),
    \end{equation*}
    so no single weighted game rejects both $L_1$ and $L_2$ following the argument for proving non-separability given in~\cite[Lemma~2]{KoberWeltge2021}. 
    Furthermore, this shows that the dimension of the game is greater than one.

    Hence, the dimension is two, but every system containing ${\cal C}_3$ needs at least three weighted games.
\end{example}

\subsection{Reduced games } \label{sec:HeuristicReducedGames}

We additionally fix the weights of smaller countries to their population data and only search for new weights for larger countries.
For the 2014 population data of the European Union, there are only ten coalitions that are winning in the (weighted) game ${\cal C}_\numSuf$, but losing in ${\cal C}_\numNec \cap {\cal P}$, see also Section~\ref{sec:results}.
We fix the weights of countries that are part of every such coalition of $\numSuf$ countries but with less than $q$ citizens.
In some sense, this approach is exactly opposite to~\cite{ChenCheungNg2019} where the population data is used as weights and only the weights for these smaller countries is changed.

Let
\[ \Wz=\{W\in\Wmin:\ |W|=\numSuf,\ p(W)<q\} \]
be the minimal winning coalitions that win only because of ${\cal C}_\numSuf$, \ie are winning but do not meet the population quota.
Furthermore, let 
\[
  \LargeCountries=\bigcup_{W\in \Wz}(N\setminus W), \qquad 
  \SmallCountries=N\setminus \LargeCountries=\bigcap_{W\in \Wz}W
\]
(with $\SmallCountries=N$ if $\Wz=\emptyset$) be a partition of the countries in smaller and larger ones (cf.~\cite{ChenCheungNg2019}). By this definition, $\LargeCountries$ contains exactly those countries that can possibly be part of an \enquote{amost-blocking minority} of three countries that together represent more than 35\% of the total population.

The sets $\LargeCountries$ and $\SmallCountries$ can be read off from sums of small populations without explicitly listing coalitions in~$\Wz$.
A country $i$ is in $\SmallCountries$ if and only if the $\numSuf$ least populous countries except $i$ have a population of at least~$q$.

\begin{definition}
    For $\SubsetLarge \subseteq \LargeCountries$, let
    \begin{align*}
        \maxPopLose_\SubsetLarge &= \max\{p(\SubsetSmall):\ \SubsetSmall\subseteq \SmallCountries,\ |\SubsetLarge|+|\SubsetSmall|\ge \numNec,\ \SubsetLarge\cup \SubsetSmall\notin\W\},\\
        \minPopWin_\SubsetLarge &= \min\{p(\SubsetSmall):\ \SubsetSmall\subseteq \SmallCountries,\ \SubsetLarge \cup \SubsetSmall\in\W \},
    \end{align*}
    with $\min\emptyset=+\infty$ and $\max\emptyset=-\infty$.
\end{definition}
By definition, we have $\maxPopLose_\SubsetLarge<\minPopWin_\SubsetLarge$ 
for all $\SubsetLarge \subseteq \LargeCountries$.

\begin{lemma}\label{lem:slice}
Let all populations be positive. For every $\SubsetLarge\subseteq \LargeCountries$ and $\SubsetSmall\subseteq \SmallCountries$ with $|\SubsetLarge|+|\SubsetSmall|\ge \numNec$,
\[
  \SubsetLarge\cup \SubsetSmall\in\W \iff p(\SubsetSmall)\ge\minPopWin_\SubsetLarge.
\]
\end{lemma}

\begin{proof}
Let $\SubsetLarge\subseteq \LargeCountries$ and $\SubsetSmall\subseteq \SmallCountries$ with $|\SubsetLarge|+|\SubsetSmall|\ge \numNec$. 
Since all populations are positive, $p(\SubsetSmall)<p(\SmallCountries)$ holds for every $\SubsetSmall\subsetneq \SmallCountries$. 
Thus, $\SubsetLarge\cup \SubsetSmall$ with $\SubsetSmall\neq \SmallCountries$ wins if and only if $p(\SubsetLarge)+p(\SubsetSmall)\ge q$. 
If $\SubsetSmall=\SmallCountries$, then $\SubsetLarge\cup \SmallCountries$ wins if and only if $p(\SubsetLarge)+p(\SmallCountries)\ge q$ or $|\SubsetLarge|+|\SmallCountries|\ge \numSuf$. 
Hence, the winning coalitions in $\{\SubsetLarge\cup \SubsetSmall: \SubsetSmall\subseteq \SmallCountries,\ |\SubsetLarge|+|\SubsetSmall|\ge \numNec\}$ are exactly those with $p(\SubsetSmall) \geq \min\set{q-p(\SubsetLarge), p(\SmallCountries)}=\minPopWin_\SubsetLarge$ if $|\SubsetLarge| + |\SmallCountries| \geq \numSuf$ and 
otherwise $p(\SubsetSmall) \geq q-p(\SubsetLarge) = \minPopWin_\SubsetLarge$.
\end{proof}

Positivity of $p$ cannot be dropped as the following example shows. 
\begin{example}
Let $p=(1,0,0)$, $q=1$, $\numNec=1$ and $\numSuf=2$.
    Then, a coalition wins if and only if it contains player $1$ or has at least two members. 
    We have $\Wz=\{\{2,3\}\}$, $\LargeCountries=\{1\}$, and $\SmallCountries=\{2,3\}$.
    For $\SubsetLarge=\emptyset$, both the winning coalition $\SmallCountries$ and the losing coalition $\{2\}$ have population $0$.
\end{example}

For simplicity, from here on, we assume that $p$ is positive.

\begin{definition}
    We call a weighted game \emph{$\LargeCountries$-reduced} if it has the form
\[
  {\cal V}=[t;\ v,\ p_\SmallCountries],\qquad v\in\mathbb{Z}^\LargeCountries_{\ge 0},\ t\in\mathbb{Z},
\]
\ie its weight on countries $i \in \SmallCountries$ are their populations $p_i$. 
\end{definition}

\begin{theorem}\label{thm:reduction}
Let all populations be positive.
Furthermore, let ${\cal V}_2,\dots,{\cal V}_k$ be $\LargeCountries$-reduced games with ${\cal V}_i=[t^{(i)}; v^{(i)}, p_\SmallCountries]$ for $i=2,\ldots,k$. Then $\W={\cal C}_\numNec\cap {\cal V}_2\cap\dots\cap {\cal V}_k$ if and only if for every $\SubsetLarge\subseteq \LargeCountries$ we have
\begin{equation}\label{eq:slices}
    \maxPopLose_\SubsetLarge < \max_{2\le i\le k} \left(t^{(i)}-v^{(i)}(\SubsetLarge)\right) \leq \minPopWin_\SubsetLarge.
\end{equation}
\end{theorem}

\begin{proof}
Coalitions with fewer than $\numNec$ members lose in $\W$ and in ${\cal C}_\numNec$. 
Let $S=\SubsetLarge\cup \SubsetSmall$ with $\SubsetLarge\subseteq \LargeCountries$, $\SubsetSmall\subseteq \SmallCountries$ and $|\SubsetLarge|+|\SubsetSmall|\ge \numNec$. 
If $S\in\W$, then $p(\SubsetSmall)\ge\minPopWin_\SubsetLarge\ge t^{(i)}-v^{(i)}(\SubsetLarge)$ by Lemma~\ref{lem:slice}, so $S$ wins in every~${\cal V}_i$. 
If $S\notin\W$, then $p(\SubsetSmall)\le\maxPopLose_\SubsetLarge<t^{(i)}-v^{(i)}(\SubsetLarge)$ for some $i$, so $S$ loses in ${\cal V}_i$. 

Conversely, if $t^{(i)}-v^{(i)}(\SubsetLarge)>\minPopWin_\SubsetLarge$, the winning coalition $\SubsetLarge\cup \SubsetSmall$ with $p(\SubsetSmall)=\minPopWin_\SubsetLarge$ loses in ${\cal V}_i$.
If $t^{(i)}-v^{(i)}(\SubsetLarge)\le\maxPopLose_\SubsetLarge$ for all $i$, the losing coalition $\SubsetLarge\cup \SubsetSmall$ with $p(\SubsetSmall)=\maxPopLose_\SubsetLarge$ wins in all ${\cal V}_i$.
\end{proof}

Since $\maxPopLose_\SubsetLarge$ and $\minPopWin_\SubsetLarge$ are integers, condition~\eqref{eq:slices} is a covering problem on at most $2^{|\LargeCountries|}$ points: 
Each $\LargeCountries$-reduced game is a function $\SubsetLarge\mapsto t-v(\SubsetLarge)$ that stays below all $\minPopWin_\SubsetLarge$ and \enquote{covers} the subsets $\SubsetLarge \subseteq \LargeCountries$ with $t-v(\SubsetLarge)\ge\maxPopLose_\SubsetLarge+1$. 
Restricting ourselves to $\LargeCountries$-reduces games to cover all subsets with  $\maxPopLose_\SubsetLarge>-\infty$, we might not be able to reach the dimension.
First, including ${\cal C}_\numNec$ in a representation may cost one game, see Example~\ref{example:fixCxCanCost}. 
Second, even when ${\cal C}_\numNec$ belongs to a minimum representation, the other games need not be $\LargeCountries$-reduced, as the following example shows.

\begin{example}\label{prop:notsharp}
Let $N=\{1,\dots,8\}$, $p=(9,9,6,3,3,2,1,1)$, $q=19$, $\numNec=4$ and $\numSuf=6$. 
We have $\Wz=\{\{3,\dots,8\}\}$, $\SmallCountries=\{3,\dots,8\}$ and $\LargeCountries=\{1,2\}$.

Then $\dim\W=2$ and $\W={\cal C}_4\cap {\cal V}$ with ${\cal V}=[18;7,7,6,3,3,2,2,2]$, so ${\cal C}_4$ belongs to a minimum representation. 
The losing coalitions $\{1,2,3\}$ and $\{1,4,5,6,7\}$ and the winning coalitions $\{1,3,4,5\}$ and $\{1,2,6,7\}$ 
have the same multiset of members, so $\W$ is not weighted and $\dim\W\ge 2$.
In the remainder of this paragraph, we argue that $\W=C_4\cap V$.
Let $a=|S\cap\{1,2\}|$ and $b=|S\cap\{7,8\}|$.
Then $v(S)=p(S)-2a+b$ holds, where $v$ denotes the weights of $\cal V$. 
It suffices to show that coalitions with at least four members satisfy $v(S)\ge 18$ if they win and $v(S)\le 17$ if they lose. 
A winning coalition with at least six members has $v(S)\ge 18$.  Let $S$ win with $|S|\ge 4$ and $p(S)\ge 19$. 
    If $a=0$, then $v(S)\ge p(S)\ge 19$. 
    If $a=1$, write $S=\{i\}\cup \SubsetSmall$ with $|\SubsetSmall|\ge 3$ and $p(\SubsetSmall)\ge 10$.
        If $b \geq 1$, then $v(S)=7+p(\SubsetSmall)+b\ge 18$. 
        Otherwise, $b = 0$ hold, so $\SubsetSmall \subseteq\set{3,4,5,6}$, so $p(\SubsetSmall)\in\{8,11,12,14\}$ and thus $p(\SubsetSmall)\ge 11$.
    In the remaining case if $a=2$, the two further members add at least $4$ to $v(\{1,2\})=14$. 
Now let $S$ lose with $|S|\in\{4,5\}$, so $p(S)\le 18$.
Then, we have $a\le 1$. 
    If $a=0$, then $v(S)\le v(\{3,\dots,8\})-2=16$. 
    If $a=1$, write $S=\{i\}\cup \SubsetSmall$ with $p(\SubsetSmall)\le 9$.
        If $b \leq 1$, then $v(S)=7+p(\SubsetSmall)+b\le 17$.
        Otherwise, we have $b=2$. then $\SubsetSmall=\{7,8\}\cup R$ with $R\subseteq\{3,4,5,6\}$, $|R|\le 2$ and $p(R)\le 7$.
        Actually, we have $p(R)\le 6$ because no such $R$ has population $7$, and thus $v(S)=7+p(R)+4\le 17$.

However, no weighted game of the form $[t;v,p_\SmallCountries]$ satisfies $\W={\cal C}_4\cap[t;v,p_\SmallCountries]$. 
For $\SubsetLarge=\emptyset$, the only winning $\SubsetSmall$ with $|\SubsetSmall|\ge 4$ is $\SmallCountries$ with a population of $16$,
and the most populous $\SubsetSmall\subsetneq \SmallCountries$ with $4$ or $5$ members has population $15=\maxPopLose_{\emptyset}$. 
For $\SubsetLarge=\{1\}$ the winning $\SubsetSmall$ with $|\SubsetSmall|\ge 3$ are those with $p(\SubsetSmall)\ge 10$ or $|\SubsetSmall|\ge 5$, 
the least populous of which have population $10$, while $\{3,6,7\}$ gives a losing coalition with $p(\SubsetSmall)=9$.
The case $\SubsetLarge=\{2\}$ is identical as $1$ and $2$ have the same population.
For $\SubsetLarge=\{1,2\}$ every $\SubsetSmall \subseteq \SmallCountries$ with $|\SubsetSmall|\ge 2$ wins, the least populous being $\{7,8\}$. 
Thus, in total, we have
\[
\begin{array}{lcccc}
\SubsetLarge & \emptyset & \{1\} & \{2\} & \{1,2\}\\
\maxPopLose_\SubsetLarge & 15 & 9 & 9 & -\infty\\
\minPopWin_\SubsetLarge & 16 & 10 & 10 & 2.
\end{array}
\]

Using Theorem~\ref{thm:reduction}, 
twice the upper and twice the lower bound,
for a single $\LargeCountries$-reduced game we get the contradiction
\[
  \minPopWin_{\{1,2\}} \geq t-(v_1 + v_2)=t-v_1+t-v_2-t> \maxPopLose_{\{1\}} + \maxPopLose_{\{2\}} - \minPopWin_{\emptyset} = 9+9-16=2=\minPopWin_{\{1,2\}}.
\]
\end{example}

\subsection{Integer weights}

Requiring integer weights and quotas with the exact populations on $\SmallCountries$, as in
the definition of $\LargeCountries$-reduced games, is a further restriction.
\begin{example} 
For $p=(14,12,7,6,5,3,3,1)$, $q=28$, $\numNec=4$ and $\numSuf=6$ we have $\LargeCountries=\{1,2\}$.
Enumerating the subsets of $\LargeCountries$, we evaluate \[
\begin{array}{lcccc}
\SubsetLarge & \emptyset & \{1\} & \{2\} & \{1,2\}\\
\maxPopLose_\SubsetLarge & 24 & 13 & 15 & -\infty\\
\minPopWin_\SubsetLarge & 25 & 14 & 16 & 4.
\end{array}
\]
We claim that the dimension of this game is at most two, but it cannot be represented as the intersection of ${\cal C}_\numNec$ and a weighted game \emph{with integer weights}.
Considering the values of $\maxPopLose_\SubsetLarge$ and $\minPopWin_\SubsetLarge$ for $\SubsetLarge \in \{\emptyset, \{1\}, \{2\}\}$, Theorem~\ref{thm:reduction} forces $t=25$, $v_1=11$ and $v_2=9$ for a single integral $\LargeCountries$-reduced game.
We get a contradiction for the bound $\minPopWin_{\{1,2\}}=4 \ngeq 5 =25-11-9= t-v_1-v_2$. In fact, the coalition $\{1,2\} \cup \{7,8\}$ loses in ${\cal C}_\numNec \cap {\cal V}_2$, but should be winning.

Instead, the fractional weights $v_1=11.5$ and $v_2=9.5$ satisfy all conditions.
\end{example}

Let $\subsetsLargeCtoCover=\{\SubsetLarge\subseteq \LargeCountries:\ \maxPopLose_\SubsetLarge>-\infty\}$ be those subsets of large countries that must be
covered and $\subsetsLargeWithFinitAlpha=\{\SubsetLarge\subseteq \LargeCountries:\ \minPopWin_\SubsetLarge<\infty\}$ those that impose
an upper bound in Theorem~\ref{thm:reduction}.
It suffices to consider subsets of these two sets based on dominance in the following sense.
If $\SubsetLarge\subseteq \SubsetLarge'$, $\maxPopLose_\SubsetLarge\le\maxPopLose_{\SubsetLarge'}$ and $t-v(\SubsetLarge')>\maxPopLose_{\SubsetLarge'}$, then $t-v(\SubsetLarge)>\maxPopLose_\SubsetLarge$ for non-negative weights $v$.
Let $\subsetsLargeCtoCoverDom$ be the non-dominated sets in $\subsetsLargeCtoCover$, keeping the larger set of countries in case of a tie.
Similarly, if $\SubsetLarge\subseteq \SubsetLarge'$, $\minPopWin_\SubsetLarge\le\minPopWin_{\SubsetLarge'}$ and $t-v(\SubsetLarge)\le\minPopWin_\SubsetLarge$, then $t-v(\SubsetLarge')\le\minPopWin_{\SubsetLarge'}$ for non-negative $v$.
Let $\subsetsLargeWithFinitAlphaDom$ be the non-dominated sets in $\subsetsLargeWithFinitAlpha$, keeping the smaller sets in case of a tie.

For $\mathcal{S}\subseteq\subsetsLargeCtoCoverDom$, a $\LargeCountries$-reduced game with integral weights and quota covering all subsets $\SubsetLarge \in \mathcal{S}$ exists if and only if the set 
\begin{equation}\label{eq:ip}
    \left\{ (t, v_\LargeCountries) \in \Z \times \Z_{\geq 0}^\LargeCountries \mid ~
        t-v(\SubsetLarge) \le\minPopWin_\SubsetLarge \, (\SubsetLarge\in\subsetsLargeWithFinitAlphaDom),\quad
        t-v(\SubsetLarge) \ge\maxPopLose_\SubsetLarge+1 \, (\SubsetLarge\in\mathcal{S})
    \right\}
\end{equation}
is non-empty.
We search for $t, v_\LargeCountries$ in~\eqref{eq:ip} with an integer programming solver, minimising $\sum_j v_j$.

Finally, we compute an optimal colouring of $\subsetsLargeCtoCoverDom$ using a \bnb scheme.
First, for all ${\cal S} \subseteq \subsetsLargeCtoCoverDom$ with $|{\cal S}|=2$, we try whether a $\LargeCountries$-reduced game can cover $\cal S$.
Next, based on a descending order of the vertices by degree in this simple graph, we search for a clique.
The size of the clique is a lower bound for the colouring. 
Afterwards, we search with a \bnb scheme for an optimal colouring of the hypergraph with vertices $\subsetsLargeCtoCoverDom$ and (hyper-)edges connecting those vertices that cannot be covered by a $\LargeCountries$-reduced game.
Note that this is the analogue construction to finding a clique in the (hyper-)graph of losing coalitions used for the lower bounds~\cite{KurzNapel2016, KoberWeltge2021}.

\section{Computational Results}
\label{sec:results}

We implemented our heuristic in Python.
The IPs are solved with Highs via scipy. We use an AMD EPYC 7742 64-core CPU. 

The voting rule of the EU Council for 2014 has the structure of \eqref{eq:rule} with $\numNec=16$, $\numSuf=25$ and $q=\lceil 0.65\,p(N)\rceil=329\,820\,795$ and populations $p$ as given in Table~\ref{tab:pop2014}. 
We have
\[
  |\Wmin| = 8\,248\,125,\qquad 
  |\Lmax| = 7\,179\,233,\qquad
  |\LmaxX| = 384\,338
\]
which is significantly less than considering all $2^{28}$ coalitions.
Note that the sizes of $\Wmin$ and $\Lmax$ are also stated in \cite{KurzNapel2016}. 

$\Wz$ consists of the complements of the ten triples of the six most populous countries with more
than $35\%$ of the population,
\[
\begin{gathered}
\{\mathrm{DE},\mathrm{FR},\mathrm{UK}\},\ 
\{\mathrm{DE},\mathrm{FR},\mathrm{IT}\},\
\{\mathrm{DE},\mathrm{FR},\mathrm{ES}\},\ 
\{\mathrm{DE},\mathrm{FR},\mathrm{PL}\},\
\{\mathrm{DE},\mathrm{UK},\mathrm{IT}\},\\
\{\mathrm{DE},\mathrm{UK},\mathrm{ES}\},\ 
\{\mathrm{DE},\mathrm{UK},\mathrm{PL}\},\
\{\mathrm{DE},\mathrm{IT},\mathrm{ES}\},\ 
\{\mathrm{DE},\mathrm{IT},\mathrm{PL}\},\
\{\mathrm{FR},\mathrm{UK},\mathrm{IT}\},
\end{gathered}
\]
so $\LargeCountries=\{\mathrm{DE},\mathrm{FR},\mathrm{UK},\mathrm{IT},\mathrm{ES},\mathrm{PL}\}$
and $\SmallCountries$ consists of the remaining $22$ countries, see also \cite{ChenCheungNg2019}. 
Of the $2^6=64$ subsets of $\LargeCountries$, for the $22$ subsets $\SubsetLarge \subseteq \LargeCountries$ with $|\SubsetLarge| \leq 2$, we have $\SubsetLarge \cup \SmallCountries \notin \W$ though they contain more than $16$ countries, \ie $|\subsetsLargeWithFinitAlpha|=42$. Furthermore, three subsets of $\LargeCountries$ contain no such losing coalition, so $61$ subsets have to be covered; 
$54$ of them remain after the dominance reduction, \ie $|\subsetsLargeCtoCoverDom| = 54$. The gaps $\minPopWin_\SubsetLarge-\maxPopLose_\SubsetLarge$ can be as small as $3$ citizens. Exact colouring via a \bnb scheme on the obtained hypergraph with $54$ vertices and $344$ edges connecting two of these takes less than ten seconds.
In particular, there is a clique of seven vertices already in the simple graph that is found by a greedy method based on the vertex degrees.

The result is a colouring with seven $\LargeCountries$-reduced games, see Tables~\ref{tab:games2014} and~\ref{tab:partition2014}, and the game~${\cal C}_\numNec$.
Note that most weights are either zero or within a similar size of the according country's population.
Thus, there are eight weighted games such that the intersection are the winning coalitions in the EU Council with the 2014 data. 
Together with the game requiring at least $\numNec=16$ countries to vote in favour, we reach the lower bound of Kober and Weltge \cite{KoberWeltge2021}.
This solves the challenge of \cite{KurzNapel2016} to determine the dimension of the voting rules of the EU Council for 2014.

\paragraph{Data availability.} 
The weights of Theorem~\ref{thm:main}, given in Table~\ref{tab:games2014} in machine-readable form and a stand-alone verification program (written without the use of generative AI) are available \cite{VerificationCode}.

\paragraph{Use of generative AI.}
For the preparation of this manuscript, we used Claude (model Opus~5.5) for writing first versions of code implementing the heuristic
and parts of this manuscript, as well as for finding the small examples in Section~\ref{sec:heuristics}. 
All parts were carefully reviewed and edited by the authors.
The authors take full responsibility for the content of this article.

\bibliographystyle{alphaurl}

\appendix

\begin{table}[ht]
    \centering
    \fitsize{
    \begin{tabular}{lr@{\hspace{2em}}lr@{\hspace{2em}}lr}
    \toprule
    Country & Population & Country & Population & Country & Population \\
    \cmidrule{1-6}
Germany (DE) & 80\,780\,000 & Czechia (CZ) & 10\,512\,419 & Croatia (HR) & 4\,246\,700 \\
    France (FR) & 65\,856\,609 & Portugal (PT) & 10\,427\,301 & Lithuania (LT) & 2\,943\,472 \\
    United Kingdom (UK) & 64\,308\,261 & Hungary (HU) & 9\,879\,000 & Slovenia (SI) & 2\,061\,085 \\
    Italy (IT) & 60\,782\,668 & Sweden (SE) & 9\,644\,864 & Latvia (LV) & 2\,001\,468 \\
    Spain (ES) & 46\,507\,760 & Austria (AT) & 8\,507\,786 & Estonia (EE) & 1\,315\,819 \\
    Poland (PL) & 38\,495\,659 & Bulgaria (BG) & 7\,245\,677 & Cyprus (CY) & 858\,000 \\
    Romania (RO) & 19\,942\,642 & Denmark (DK) & 5\,627\,235 & Luxembourg (LU) & 549\,680 \\
    Netherlands (NL) & 16\,829\,289 & Finland (FI) & 5\,451\,270 & Malta (MT) & 425\,384 \\
    Belgium (BE) & 11\,203\,992 & Slovakia (SK) & 5\,415\,949 &  &  \\
    Greece (GR) & 10\,992\,589 & Ireland (IE) & 4\,604\,029 &  &  \\
    \bottomrule
  \end{tabular}}
    \caption{Population data ($2014$) as in \cite{KurzNapel2016,ChenCheungNg2019}. }
    \label{tab:pop2014}
\end{table}

\begin{table}[ht]
  \centering
  \fitsize{
  \begin{tabular}{c|r|rrrrrr} \toprule
    & quota & Germany & France & United Kingdom & Italy & Spain & Poland \\ \midrule
$\V_2$ & 241\,861\,457 & 80\,770\,349 & 65\,856\,530 & 30\,533\,730 & 30\,533\,730 & 30\,533\,730 & 30\,533\,730 \\ $\V_3$ & 221\,475\,661 & 80\,757\,531 & 0 & 35\,607\,697 & 60\,782\,563 & 35\,607\,697 & 35\,607\,697 \\ $\V_4$ & 182\,342\,377 & 80\,775\,777 & 0 & 64\,308\,150 & 0 & 32\,082\,110 & 32\,082\,110 \\ $\V_5$ & 277\,972\,216 & 80\,778\,980 & 46\,507\,585 & 46\,083\,986 & 46\,507\,585 & 46\,507\,585 & 38\,495\,468 \\ $\V_6$ & 235\,263\,586 & 0 & 52\,079\,503 & 64\,308\,206 & 60\,782\,604 & 46\,507\,710 & 38\,495\,609 \\ $\V_7$ & 184\,732\,063 & 0 & 65\,856\,385 & 0 & 60\,782\,405 & 46\,507\,550 & 38\,495\,410 \\ $\V_8$ & 188\,197\,102 & 0 & 65\,852\,334 & 64\,299\,121 & 0 & 46\,455\,873 & 38\,486\,515 \\ \bottomrule
  \end{tabular}
  }
  \caption{System of $8$ games for 2014 EU population data: $\V_1=[16;1, \dots, 1]$ and $\V_2, \dots, \V_8$ with quotas and weights for the six largest countries forming $\LargeCountries$.
  Weights of remaining countries are their populations, see Table~\ref{tab:pop2014}.} \label{tab:games2014}
\end{table}

\begin{table}[ht]
    \centering
    \fitsize{
    \begin{tabular}{r|l}
        $\V_2$ & $\{\set{\mathrm{\bf UK, IT, ES, PL}}, \set{\mathrm{UK, IT, ES}}, \set{\mathrm{FR, UK, IT, ES, PL}}, \set{\mathrm{UK, IT, PL}}, \set{\mathrm{DE, UK, IT, ES, PL}}, $ 
            \\ &
            $\phantom{\{} \set{\mathrm{UK, ES, PL}}, \set{\mathrm{IT, ES, PL}}, \set{\mathrm{UK, IT}}, \set{\mathrm{UK, ES}}, \set{\mathrm{IT, ES}}, \set{\mathrm{UK, PL}}, \set{\mathrm{IT, PL}}, \set{\mathrm{ES, PL}}\}$ \\
        $\V_3$ & $\{\set{\mathrm{\bf FR, UK, ES, PL}}, \set{\mathrm{FR, UK, ES}}, \set{\mathrm{FR, UK, PL}}, \set{\mathrm{FR, ES, PL}}, \set{\mathrm{DE, FR, UK, ES, PL}}, $ 
            \\ &
            $\phantom{\{} \set{\mathrm{FR, UK}}, \set{\mathrm{FR, IT}}, \set{\mathrm{FR, ES}}, \set{\mathrm{FR, PL}}\}$\\
        $\V_4$ & $\set{\set{\mathrm{\bf FR, IT, ES, PL}}, \set{\mathrm{FR, IT, ES}}, \set{\mathrm{FR, IT, PL}}, \set{\mathrm{DE, FR, IT, ES, PL}}} $\\
        $\V_5$ & $\{\set{\mathrm{\bf FR, UK, IT, PL}}, \set{\mathrm{FR, UK, IT}}, \set{\mathrm{FR, UK, IT, ES}}, \set{\mathrm{DE, FR}}, \set{\mathrm{DE, UK}}, \set{\mathrm{DE, IT}},  $
            \\ &
            $\phantom{\{}\set{\mathrm{DE, FR, UK, IT}}, \set{\mathrm{DE, ES}}, \set{\mathrm{DE, PL}} \}$\\
        $\V_6$ & $\{\set{\mathrm{\bf DE, FR, ES, PL}}, \set{\mathrm{DE, FR, ES}}, \set{\mathrm{DE, FR, PL}}, \set{\mathrm{DE, FR, IT, PL}}, \set{\mathrm{DE, ES, PL}},$
            \\ &
            $\phantom{\{} \set{\mathrm{DE, FR, UK, PL}}, \set{\mathrm{DE, FR, IT, ES}}, \set{\mathrm{DE, FR, IT}}, \set{\mathrm{DE, FR, UK, ES}}, \set{\mathrm{DE, FR, UK}} 
            \}$ \\
        $\V_7$ & $\set{\set{\mathrm{\bf DE, UK, ES, PL}}, \set{\mathrm{DE, UK, ES}}, \set{\mathrm{DE, UK, PL}}, \set{\mathrm{DE, UK, IT, PL}}, \set{\mathrm{DE, UK, IT}}, \set{\mathrm{DE, UK, IT, ES}}}$\\
        $\V_8$ & $\set{\set{\mathrm{\bf DE, IT, ES, PL}}, \set{\mathrm{DE, IT, ES}}, \set{\mathrm{DE, IT, PL}}}$ \\

    \end{tabular}}
    \caption{Partition of $\subsetsLargeCtoCoverDom$ for the EU data of 2014. The first coalition of every set (in bold) form a clique in the colouring graph.}
    \label{tab:partition2014}
\end{table}

\end{document}